\documentclass[12pt]{amsart}
\usepackage{amsmath}

\newtheorem{theo}{Theorem}[section]
\newtheorem{prop}[theo]{Proposition}
\newtheorem{lemma}[theo]{Lemma}
\newtheorem{cor}[theo]{Corollary}
\newtheorem{observation}[theo]{Observation}

\title[Multiparameter ergodic averages for amenable groups]{Multiparameter ergodic averages for two commuting actions of an amenable group}

\author{JOHN T. GRIESMER}

\address{Department of Mathematics, The Ohio State University, 231 W. 18th Ave., Columbus, OH 43212}
\email{griesmer@math.osu.edu}

\subjclass[2000]{37A05}

\date{\today}

\begin{document}

\begin{abstract}
We find limits of some multiple ergodic averages, generalizing a result of Bergelson to the setting of two commuting transformations and actions of amenable groups.
\end{abstract}

\maketitle

\section{Introduction}
In order to give an ergodic-theoretic proof of Szemer{\'e}di's theorem on arithmetic progressions, Furstenberg (\cite{F1}) analyzed the behavior of the averages
\begin{align}\label{1}
\frac{1}{N-M}\sum_{n=M}^{N-1} \int f \cdot f\circ T^n\cdot f\circ T^{2n}\cdot \cdots \cdot f\circ T^{kn}\, d\mu,
\end{align}
where $T:X\to X$ is a transformation of a measure space $(X,\mathcal B,\mu)$ which preserves a probability measure $\mu,$ and $f\in L^\infty(\mu).$  This raised the problem of deciding whether such averages converge as $N-M\to \infty,$ and if they do converge, identifying the limit.   The case $k=2$ was settled affirmatively in \cite{F1}, and the case $k=3$ was settled affirmatively by Conze and Lesigne in \cite{CL84} and \cite{CL88}.  Zhang  \cite{QZ} generalized the results of \cite{CL88} to the case where $T^n,T^{2n},$ and $T^{3n}$ are replaced by $T^n,S^n,$ and $R^n,$ where $R,S,$ and $T$ are commuting, ergodic transformations, and the transformations $R^{-1}S, S^{-1}T,$ and $T^{-1}R$ are assumed to be ergodic.  Finally, the question of convergence in (\ref{1}) was settled affirmatively for all $k$ by Host and Kra in \cite{HK}, and independently by Ziegler in \cite{Z}.  Each of the aforementioned papers uses a method of reducing the problem to studying the appropriate \textit{characteristic factors}, a notion which we explain in the next section.   In these proofs, identifying the appropriate characteristic factors is an essential step in identifying the limit.  Recently, Tao has generalized this convergence result to the case where the powers of $T$ in (\ref{1}) are replaced by commuting transformations $T_1,\dots, T_k.$  Interestingly, Tao's proof (\cite{T}) does not use the method of characteristic factors, and does not identify the limit of the averages.  T. Austin (\cite{A}) has found another, more general proof of Tao's theorem which does not identify the characteristic factors.  The proof in \cite{A} uses properties of diagonal measures, as in the proof of Theorem 2.3 and Corollary 2.4 of \cite{F1}.

The following related result is due to Bergelson, and appears as Theorem 5.3 in \cite{B}.  Here a \textit{syndetic} subset of $\mathbb Z^2$ is one that meets every square of side length $L$ for some $L>0.$
\begin{theo}\label{B1}  Let $(X,\mathcal B,\mu, T)$ be an invertible probability measure preserving system.  Then:
\begin{enumerate}
\item For any $f_1,f_2,f_3\in L^\infty(X,\mathcal B,\mu),$
\begin{align*}
\lim_{N-M\to \infty} \frac{1}{(N-M)^2} \sum_{n,m=M}^{N-1} f_1(T^n x)f_2(T^mx)f_3(T^{n+m}x)
\end{align*}
exists in $L^2(\mu).$
\item  For any $A\in \mathcal B$ with $\mu(A)>0$
\begin{align*}
\lim_{N-M\to \infty} \frac{1}{(N-M)^2}\sum_{n,m=M}^{N-1} \mu(A\cap T^nA\cap T^m A\cap T^{n+m}A)\geq \mu(A)^4.
\end{align*}
\end{enumerate}
\end{theo}
From this fact, Bergelson derives the following combinatorial corollary which generalizes Khintchine's refinement of the Poincar{\'e} recurrence theorem.
\begin{cor}\label{Bcor1}  For $(X,\mathcal B,\mu ,T)$ as in Theorem \ref{B1}, for any $0<\lambda<1$ and any $A\in \mathcal B$ with $\mu(A)>0,$ the set
\begin{align*}
\{(n,m)\in \mathbb Z^2: \mu\left((A\cap T^nA)\cap T^m(A\cap T^nA)\right)>\lambda\mu(A)^4\}
\end{align*}
is syndetic.
\end{cor}

Using Furstenberg's correspondence principle (see \cite{Bert}) 
one can derive a combinatorial consequence of Corollary \ref{Bcor1}, which says roughly that given a set of positive density $\delta$ in $\mathbb Z,$ there is a syndetic set of $(n,m)\in \mathbb Z^2$ so that for all such $(n,m),$ configurations of the form $(a,a+n,a+m,a+n+m)$ appear as frequently as you would expect in a random set having density $\delta.$  The precise statement follows.

\begin{theo}\label{random}  Let $A\subseteq \mathbb Z$ such that $d^*(A)=\limsup_{N-M\to \infty} \frac{|A\cap [M,N-1]|}{N-M}=\delta.$  Then for all $\varepsilon>0,$ the set of $(n,m)$ such that
$$
d^*(A\cap (A-n)\cap (A-m)\cap (A-n-m))> \delta^4-\varepsilon
$$
is syndetic.
\end{theo}
It is natural to conjecture that a similar result holds for sets of positive density in $\mathbb Z^2,$ and that such a result may be proved by establishing convergence for averages of the form 
\begin{align}\label{avgs}
\frac{1}{(N_1-M_1)(N_2-M_2)}\sum_{n=M_1}^{N_1-1}\sum_{m=M_2}^{N_2-1} f_1(x) f_2(T^n x) f_3(S^m x) f_4(T^nS^m x),
\end{align} where $T$ and $S$ are commuting, measure preserving transformations of a probability measure space.  In this paper we establish convergence of the averages in (\ref{avgs}) and use this to deduce a combinatorial consequence for sets of positive density in $\mathbf Z^2.$

We actually consider a more general situation, which may be motivated by studying configurations in sets of positive density in an amenable group $G,$ and in the cartesian square $G\times G.$  We say a group is \textit{amenable} if there is a sequence of finite sets $\Phi_n \subseteq G$ such that for all $g\in G, \lim_{n\to \infty} \frac{|\Phi_n\cap g\Phi_n|}{|\Phi_n|}=1.$  Such a sequence is called a \textit{left F{\o}lner sequence}, and these sequences allow one to define a translation-invariant notion of density in the group $G.$  Analogously, we say that $\Phi_n$ is a \textit{right} F{\o}lner sequence if for all $g\in G, \lim_{n\to \infty} \frac{|\Phi_n\cap \Phi_ng|}{|\Phi_n|}=1,$ and we say that $\Phi_n$ is a \textit{two-sided} F{\o}lner sequence if it is both a left- and a right F{\o}lner sequence.

Supposing that $G$ is an amenable group and that $g\mapsto S_g,g\mapsto T_g$ are commuting actions of $G$ by measure preserving transformations on a standard probability space $(X,\mathcal B,\mu),$ we demonstrate that the limit
\begin{align*}
\lim_{n\to \infty} \frac{1}{|\Phi_n||\Psi_n|}\sum_{g\in \Phi_n, h\in \Psi_n} f_1(T_gx)f_2(S_hx)f_3(T_gS_h x)
\end{align*}
exists for all two-sided F{\o}lner sequences $\Phi,\Psi$ in $G$ and that the limit is independent of the choice of F{\o}lner sequences.  We identify the limit, and also give a combinatorial inequality analogous to Corollary \ref{Bcor1}.  Further, following \cite{BR}, we give necessary and sufficient conditions for the limit to be equal to a constant almost everywhere.

Our main result is the following theorem.

\begin{theo}\label{main}  Let $G$ be a countable amenable group, let $\mathbf X=(X,\mathcal B, \mu)$ be a standard probability space with  probability measure $\mu,$ and commuting measure preserving actions $T,S$ of $G$ on $X.$  Then

(1) For all $f_1,f_2,f_3\in L^\infty (\mu),$ and all two-sided F{\o}lner sequences $\Phi,\Psi$ in $G,$ the limit
\begin{align}\label{lim}
L=\lim_{n\to \infty} \frac{1}{|\Phi_n||\Psi_n|} \sum_{(g,h)\in \Phi_n \times \Psi_n} f_1(T_g x) f_2(S_hx) f_3(T_gS_hx)
\end{align}
exists in $L^2(\mu).$  

(2) $L$ is equal to a constant $\mu$-almost everywhere for all $f_1,f_2,f_3\in L^\infty(\mu)$ if and only if $T\times T$ and $S\times S$ are ergodic. 
\end{theo}
The condition that $T\times T$ is ergodic is of course equivalent to the condition that $T$ is weakly mixing.

It is not much of a sacrifice to restrict our attention to two-sided F{\o}lner sequences, for every countable amenable group admits two-sied F{\o}lner sequences - see \cite{G} for a proof.

Theorem \ref{main} has a corollary analogous to Corollary \ref{Bcor1}.

\begin{cor}\label{large}  Let $(X,\mathcal B,\mu),T$ and $S$ be as in Theorem \ref{main}.  Suppose $f\in L^\infty(\mu)$ is a nonnegative function.  Then  for all two-sided F{\o}lner sequences $\Phi,\Psi$ 
\begin{align}\label{large1}
\lim_{n\to \infty} \frac{1}{|\Phi_n||\Psi_n|} \sum_{(g,h)\in \Phi_n\times \Psi_n} \int fT_gfS_hfT_gS_hf\, d\mu\geq \left(\int f\, d\mu \right)^4.
\end{align}
\end{cor}

We derive several consequences from Corollary \ref{large}.  First there is the natural combinatorial result about sets of positive density in $G\times G$.  In a general group $G$, we say a set $S\subseteq G$ is \textit{left syndetic} if there is a finite set $F\subseteq G$ such that $G=\bigcup_{g\in F} gS,$ and we say that $S$ is \textit{right} syndetic if there is a finite set $F$ such that $G=\bigcup_{g\in F} Sg.$

\begin{cor}\label{combinatorial} Let $E\subseteq G\times G,$ and let $\Phi$ be a F{\o}lner sequence in $G\times G$ with $\limsup_{n\to \infty} \frac{|E\cap \Phi_n|}{|\Phi_n|}=\delta>0.$  Then for all $\varepsilon>0,$ the set
\begin{align*}
\left\{(g,h): \limsup_{n\to \infty} \frac{|E\cap E(g,1_G)\cap E(1_G,h)\cap E(g,h)\cap\Phi_n|}{|\Phi_n|}>\delta^4-\varepsilon\right\}
\end{align*}
is both left- and right syndetic in $G\times G$
\end{cor}
A combinatorial corollary about partitions of subsets in $G\times G\times G$ will also be derived in Section \ref{Combinatorial Corollaries}.

In Section \ref{example} we specialize to the case $G=\mathbb Z,$ and use the main theorem to construct a class of examples which illustrates the difference between the case of powers of a single transformation as studied in \cite{HK} and the case of commuting transformations.

\section{Acknowledgements.}  The author thanks Vitaly Bergelson for posing the problem of establishing Theorem \ref{main} and for helpful comments and advice.  Thanks are due Alexander Leibman for suggesting the example in Section \ref{example}.

\section{Notation and background.}  Throughout the paper, we fix a probability space $(X,\mathcal B,\mu).$  We may assume, without loss of generality, that $(X,\mathcal B,\mu)$ is a Lebesgue space - that is, measure-theoretically isomorphic to the unit interval with Lebesgue measure and the completed Borel $\sigma$-algebra, together with possibly countably many atoms.  (See \cite{Roy}, Chapter 15, Theorem 4, and Chapter 15, Theorem 20, for justification.)  In particular we can assume that $X$ is compact metric.

As mentioned above, every discrete amenable group has a two-sided F{\o}lner sequence.  We write $\Phi$ for a F{\o}lner sequence $\{\Phi_n\}_{n\in \mathbb N}.$  If $\Phi$ is a F{\o}lner sequence and $g\in G,$ we write $g\Phi$ and $\Phi g$ for the sequences $\{g\Phi_n\}, \{\Phi_n g\}$ respectively, and note without proof that if $\Phi$ is a (left-, right-, two-sided) F{\o}lner sequence then $g\Phi$ and $\Phi g$ are both (left-, right-, two-sided) F{\o}lner sequences.

We denote by $T$ or $S$ an action of an amenable discrete group $G$ on $(X,\mathcal B,\mu),$ and write $T_g$ for the transformation corresponding to $g,$ so that $T_gT_h=T_{gh}$ for all $g,h\in G.$  Note that this induces an anti-action of $G$ on $L^2(\mu)$ by $T_gf=f\circ T_g$ for $f\in L^2(\mu).$  We will use $T_gf$ to denote $f\circ T_g.$   Write $T\times T$ for the action of $G$ on $X\times X$ given by $g\mapsto T_g\times T_g.$  

We briefly review the notions of a factor of a dynamical system, conditional expectation, and disintegration of measures.  Write $\mathbf X=(X,\mathcal B,\mu, T,S)$ for a dynamical system where $(X,\mathcal B,\mu)$ is as above, and $T$ and $S$ are commuting actions of  $G$ on $(X,\mathcal B,\mu)$ by measure preserving transformations of $(X,\mathcal B,\mu).$  We say that $\mathbf Y=(Y,\mathcal D,\nu, T',S')$ is a \textit{factor} of $\mathbf X$ if there exists $\pi:X\to Y$ with $\pi^{-1}(\mathcal D)\subseteq \mathcal B, \mu(\pi^{-1}(A))=\nu(A)$ for all $A\in \mathcal D,$ and $T_g'\pi(x)=\pi(T_gx)$ and $S_g'\pi(x)=\pi(S_gx)$ for $\mu$-almost all $x.$   One can identify $\mathcal D$ with the $\sigma$-algebra $\pi^{-1}(\mathcal D).$

Given a $T$- and $S$-invariant, countably generated sub $\sigma$-algebra $\mathcal B'\subset \mathcal B,$ there is always be a corresponding factor $(Y,\mathcal D,\nu,T',S')$ so that $\mathcal B'=\{\pi^{-1}(D):D\in \mathcal D\}$ (up to sets of measure $0$), and we will abuse notation and identify $\mathcal B'$ and $\mathcal D.$

If $\mathcal A$ is a sub $\sigma$-algebra of $\mathcal B,$ and $f\in L^2(\mu),$ one can define the  conditional expectation (relative to $\mu$) $\mathbb E(f|\mathcal A)$ by $\mathbb E(f|\mathcal A)(x):=Pf(x),$ where $Pf$ is the orthogonal projection (in $L^2(\mu)$) of $f$ onto the subspace of $\mathcal A$-measurable functions.  (This can be extended to a map on $L^1(\mu),$ but we have no need for this here.)

The conditional expectation can be used to define disintegration of measures, as in \cite{F2}, Theorem 5.8.  If $\mathcal A$ is a countably generated sub $\sigma$-algebra of $\mathcal B,$  then there is a measurable map $x\mapsto \mu_x$ from $X$ to the space of regular Borel measures on $X$ (with the $\text{weak}^*$ topology), such that
\begin{enumerate}
\item[i] For every $f\in L^1(\mu), f\in L^1(X,\mathcal B,\mu_x)$ for $\mu$-almost every $x,$ and $\mathbb E(f|\mathcal A)(x)=\int f\, d\mu_x$ for $\mu$-almost every $x.$
\item[ii] $\int \int f \, d\mu_x \, d\mu(x)= \int f\, d\mu$ for every $f\in L^1(X,\mathcal B,\mu).$
\end{enumerate}
If $(Y,\mathcal D,\nu,T',S')$ is a factor of $\mathbf X$ with factor map $\pi,$ we use $x\mapsto \mu_{\pi(x)}$ to denote the disintegration of $\mu$ over $\pi^{-1}(\mathcal D).$  Note that $\int f d\mu_{\pi(T_gx)}=\int f\circ T_g \, d\mu_{\pi(x)}$ for $\mu$-almost every $x.$

Let $\mu=\int \mu_{\pi(x)} \, d\mu(x)$ be the disintegration of $\mu$ over $\mathcal D,$ so that $x\mapsto \mu_{\pi(x)}$ is a $\mathcal D$-measurable function, and $\int f\, d\mu_x= \mathbb E(f|\mathcal D)$ for $\mu$-almost every $x.$  One defines the \textit{relative product} measure $\mu \underset{\mathcal D}\times \mu$ on $(X\times X, \mathcal B\otimes \mathcal B)$ by  $\mu \underset{\mathcal D}\times \mu=\int \mu_{\pi(x)} \times \mu_{\pi(x)} \, d\mu(x).$  This means that $\int f \,d \mu \underset{\mathcal D}\times \mu=\int \int f \, d\mu_{\pi(x)}\times \mu_{\pi(x)}\, d\mu(x)$ for bounded $\mathcal B\otimes \mathcal B$-measurable $f.$  Equivalently, one may define $\mu \underset{\mathcal D}\times \mu$ by $\int f\otimes g d\mu \underset{\mathcal D}\times \mu=\int \mathbb E(f|\mathcal D)\mathbb E(g|\mathcal D)\, d\mu$ for bounded $f,g.$  The relative product system $\mathbf X\underset{\mathbf Y}\times \mathbf X$ is defined as the system $(X\times X, \mathcal B\otimes \mathcal B,\mu \underset{\mathcal D}\times \mu, T\times T,S\times S).$

If $R$ is a measure preserving action of $G$ on $(X,\mathcal B,\mu),$ define $\mathcal I_R$ to be the $\sigma$-algebra of $R$-invariant sets (up to sets of measure $0$): $\mathcal I_R=\{A\in \mathcal B:\mu(R_gA\triangle A)=0 \text{ for all } g\in G\}.$   Let $\mathbf I_R=(I_R,\mathcal I_R,\mu, G)$ be the factor of $\mathbf X$ corresponding to $\mathcal I_R.$  We will mainly be concerned with $\mathbf I_T$ and $\mathbf I_S.$  Let $\tau:X\to I_T$ and $\sigma:X\to I_S$ denote the  respective factor maps.

We will use a version of the ergodic theorem for amenable groups.

\begin{theo}\label{ergodic}  Let $\mathbf X=(X,\mathcal B,\mu, T)$ be a measure preserving system with $T$ a measure preserving action of $G$ on $X.$  Then for all two-sided F{\o}lner sequences  $\Phi$ in $G,$ and all $f\in L^2(\mu),$ the expectation $\mathbb E(f|\mathcal I_T)$ is given by
\begin{align*}
\mathbb E(f|\mathcal I_T)=\lim_{n\to \infty} \frac{1}{|\Phi_n|} \sum_{g\in \Phi_n} f\circ T_g=\lim_{n\to \infty} \frac{1}{|\Phi_n|} \sum_{g\in \Phi_n} f\circ T_g^{-1}
\end{align*}
where the convergence is in the sense of $L^2(\mu).$
\end{theo}

\section{Proof of the main theorem.}

We restate our main theorem for convenience.

\setcounter{section}{1}

\setcounter{theo}{3}

\begin{theo}  Let $G$ be a countable amenable group, let $\mathbf X=(X,\mathcal B, \mu)$ be a standard probability space with  probability measure $\mu,$ and commuting measure preserving actions $T,S$ of $G$ on $X.$  Then

(1) For all $f_1,f_2,f_3\in L^\infty (\mu),$ and all two-sided F{\o}lner sequences $\Phi,\Psi$ in $G,$ the limit
\begin{align}\label{lim'}
L=\lim_{n\to \infty} \frac{1}{|\Phi_n||\Psi_n|} \sum_{(g,h)\in \Phi_n \times \Psi_n} f_1(T_g x) f_2(S_hx) f_3(T_gS_hx)
\end{align}
exists in $L^2(\mu).$  

(2) $L$ is equal to a constant $\mu$-almost everywhere for all $f_1,f_2,f_3\in L^\infty(\mu)$ if and only if $T\times T$ and $S\times S$ are ergodic. 
\end{theo}

\setcounter{section}{4}

\setcounter{theo}{0}

Our proof will follow the proof of Theorem \ref{B1} found in \cite{B}, which is an example of the method of characteristic factors.  See \cite{F3},\cite{FW}, or \cite{HK} for a general description of this method.  Briefly, if $(X,\mathcal B, \mu, T)$ is a measure preserving system as in Theorem \ref{B1}, we say that a factor $\mathbf Y=(Y,\mathcal D,\nu,T')$ of a $\mathbf X$ is a \textit{characteristic factor for the scheme} $(T^n,T^m,T^{n+m})$ if for all $f_1,f_2,f_3\in L^\infty(\mu),$ 
$$
\lim_{N-M\to \infty} \frac{1}{N-M}\sum_{n=M}^{N-1} T^n f_1 S^mf_2T^{n+m}f_3-T^{n} \mathbb E(f_1|\mathcal D)T^m\mathbb E(f_2|\mathcal D)T^{n+m}\mathbb E(f_3|\mathcal D)=0.
$$
Thus, if $\mathbf Y$ is a characteristic factor for the scheme $(T^n,T^m,T^{n+m})$ one need only establish the convergence result for functions defined on $\mathbf Y.$  In \cite{B} it is shown that the Kronecker factor is a characteristic factor for the scheme $(T^n,T^m,T^{n+m})$ and convergence is established using properties of eigenfunctions.  For the sake of brevity, we define the notion of characteristic factor precisely only for our scenario.  

\textbf{Definition.}  Let $\mathbf X,$ $T$ and $S$ be as in the statement of Theorem \ref{main}.  We say that a factor $\mathbf Y$ of $\mathbf X$ is a \textit{characteristic factor for the scheme} $(T_g,S_h,T_gS_h)$ if for all functions $f_1,f_2,f_3\in L^\infty(\mu)$ and two-sided F{\o}lner sequences $\Phi,\Psi,$
\begin{align*}
\lim_{n\to \infty} \frac{1}{|\Phi_n||\Psi_n|}\sum_{(g,h)\in \Phi_n\times \Psi_n} T_gf_1S_hf_2T_gS_hf_3-T_g\mathbb E(f_1|\mathcal D)S_h \mathbb E(f_2|\mathcal D) T_gS_h\mathbb E(f_3|\mathcal D)=0
\end{align*}
\medbreak

Using similar methods to \cite{BMZ}, we will find a characteristic factor for the scheme $(T_g,S_h,T_gS_h).$    However, in the course of the proof, we only use \textit{partial} characteristic factors:  we find factors $\mathbf Y_1,\mathbf Y_2$ with corresponding invariant $\sigma$-algebras $\mathcal D_1,\mathcal D_2\subseteq \mathcal B$ so that
\begin{align*}
\lim_{n\to \infty} \frac{1}{|\Phi_n||\Psi_n|}\sum_{(g,h)\in \Phi_n\times \Psi_n} T_gf_1S_hf_2T_gS_hf_3-T_g\mathbb E(f_1|\mathcal D_1)S_h \mathbb E(f_2|\mathcal D_2) T_gS_hf_3=0.
\end{align*}
As with most applications of characteristic factors, we need a van der Corput lemma.  The following appears as Lemma 4.2 in \cite{BMZ}.
\begin{lemma}\label{vdc}  Suppose that $\{u_g:g\in G\}$ is a bounded set in a Hilbert space $H$, and that $\Phi$ is a left F{\o}lner sequence in $G$.  If
\begin{align*}
\lim_{n\to \infty} \frac{1}{|\Phi_n|^2} \left(\limsup_{m\to \infty} \frac{1}{|\Phi_m|} \sum_{g\in \Phi_m} \sum_{h,k\in \Phi_n} \left\langle u_{hg},u_{kg}\right\rangle\right)=0,
\end{align*}
then $\lim_{n\to \infty} \|\frac{1}{|\Phi_n|}\sum_{g\in \Phi_n} u_g\|=0.$
\end{lemma}

We will actually apply a corollary of Lemma \ref{vdc}.

\begin{cor}\label{vdc'}  Suppose that $\{u_{g,h}:g,h\in G\}$ is a bounded set in a Hilbert space $H$, and that $\Phi, \Psi$ are two-sided F{\o}lner sequences in $G.$  If
\begin{align*}
\lim_{m\to \infty} \frac{1}{|\Phi_m|^2|\Psi_m|^2}\left(\limsup_{m\to \infty} \frac{1}{|\Phi_m||\Psi_m|} \sum_{\substack{ j,j'\in \Phi_m\\ k,k'\in \Phi_m}}\sum_{\substack{g\in \Phi_n\\ h\in \Psi_n}} \langle u_{j'g,k'h}, u_{jg,kh}\rangle \right)=0
\end{align*}
then $\lim_{n\to \infty} \|\frac{1}{|\Phi_n||\Psi_n|} \sum_{{g,h}\in \Phi_n\times \Psi_n} u_{g,h}\|=0.$
\end{cor}

In order to describe our partial characteristic factors, we use the construction of the maximal isometric extension of a factor of $\mathbf X.$  If $\mathbf Y$ is a factor of $\mathbf X=(X,\mathcal B,\mu,T)$ let $W_\mathbf Y$ denote the closed $T$-invariant subspace of $L^2(\mu)$ spanned by functions of the form $k(x)=\int H(x,z)\phi(z)\, d\mu_{\pi(x)}(z),$ where $H\in L^\infty (\mu\underset{\mathcal D}\times \mu)$ is $T\times T$-invariant (with respect to $\mu\underset{\mathcal D}\times \mu$), and $\phi\in L^\infty(\mathbf X).$ In fact, $W_\mathbf Y$ is a closed algebra of functions, and therefore determines a factor of $\mathbf X,$ but this fact will not be needed in our proofs.  To check that $W_\mathbf Y$ is $T$-invariant, we write

\begin{align}\label{composition}
\begin{split}
\int H(Tx,z)\phi(z)\, d\mu_{\pi(Tx)}(z)&=\int H(Tx,Tz)\phi(Tz)\, d\mu_{\pi(x)}(z)\\
&= \int H(x,z) \phi(Tz)\, d\mu_{\pi(x)}(z).
\end{split}
\end{align}

The next lemma is essentially Lemma 7.6 from \cite{F2}, adapted to the setting of amenable groups.

\begin{lemma}\label{wm}  Let $\mathbf Y$ be a factor of $\mathbf X.$  If $f\in L^\infty(\mathbf X)$ and $f\perp W_\mathbf Y,$ then for all two-sided F{\o}lner sequences $\Phi,$ 
\begin{align*}
\lim_{n\to\infty} \frac{1}{|\Phi_n|} \sum_{g\in \Phi_n} \|\mathbb E(fT_gf|\mathcal D)\|=0.
\end{align*}
\end{lemma}

\textit{Proof.} Let $\Phi$ be a two-sided F{\o}lner sequence.  It suffices to show that
\begin{align*}
\lim_{n\to \infty} \frac{1}{|\Phi_n|} \sum_{g\in \Phi_n} \|\mathbb E(fT_gf|\mathcal D)\|^2=0.
\end{align*}  
Write $\|E(fT_gf|\mathcal D)\|^2=\int \int fT_gf\, d\mu_{\pi(x)}\int \bar fT_g\bar f \, d\mu_{\pi(x)} \, d\mu(x),$ so that
\begin{align}\label{wm1}
\|\mathbb E(fT_gf|\mathcal D)\|^2=\int \int f(w)f(T_gw)\bar f(T_g z)\bar f(z)\, d\mu_{\pi(x)}\times \mu_{\pi(x)} (w,z)\, d\mu(x).
\end{align}
Averaging over $\Phi$ in (\ref{wm1}) yields
\begin{align}\label{wm2}
\lim_{n\to \infty} \frac{1}{|\Phi_n|} \sum_{g\in \Phi_n} \|\mathbb E(fT_gf|\mathcal D)\|^2=\int \int f(w) H(w,z) \bar f(z)\, d\mu_{\pi(x)}\times \mu_{\pi(x)} (w,z)\, d\mu(x),
\end{align}
where $H$ is the projection of $f\otimes \bar f$ on the space of $T\times T$-invariant functions on $\mathbf X\underset{\mathbf Y}\times \mathbf X.$  Changing the order of integration, we rewrite (\ref{wm2}) as
\begin{align}\label{wm3}
\begin{split} \lim_{n\to \infty} \frac{1}{|\Phi_n|} \sum_{g\in \Phi_n} &\|\mathbb E(fT_gf|\mathcal D)\|\\
&= \int \int f(w) \int H(w,z) \bar f(z)\, d\mu_{\pi(x)}(z)\, d\mu_{\pi(x)}(w)\, d\mu(x)\\
&=\int f(x) \int H(x,z)\bar f(z)\, d\mu_{\pi(x)}(z)\, d\mu(x),
\end{split}
\end{align}
where the last equality is justified by the disintegration $\mu=\int \mu_{\pi(x)}\, d\mu(x).$  By assumption, $f$ is orthogonal to the function given by the inner integral $\int H(x,z)\bar f(z)\, d\mu_x(w),$ so the last integral in (\ref{wm3}) is $0.$ $\square$
\medbreak

\subsection{Reduction to partial characteristic factors.}  Let $W_{T/S}$ denote the closed subspace of $L^2(\mu)$ spanned by functions of the form $k(x)=\int H(x,z)\phi(z)\, d\mu_{\sigma(x)}(z),$ where $H$ is $T\times T$-invariant (with respect to $\mu \underset{\mathcal I_S}\times \mu,$) and $\phi\in L^\infty(\mu).$  Let $W_{S/T}$ denote the closed subspace of $L^2(\mu)$ spanned by functions of the form $k(x)=\int H(x,z)\phi(z)\, d\mu_{\tau(x)}(z),$ where $H$ is $S\times S$-invariant (with respect to $\mu \underset{\mathcal I_S}\times \mu$).  Let $P:L^2(\mu)\to W_{T/S}$ denote the orthogonal projection onto $W_{T/S},$ and let $Q:L^2(\mu)\to W_{S/T}$ denote the projection onto $W_{S/T}.$

The spaces $W_{T/S}$ and $W_{S/T}$ will serve as partial characteristic factors, in that we will reduce the problem of showing that the limit exists in (\ref{lim}) to the case where $f_1\in W_{T/S}$ and $f_2\in W_{S/T}.$  We will not, however, show that these spaces are closed under pointwise multiplication, so we are not considering them as factors at present.  Once convergence in (\ref{lim}) is established, we will see that $W_{T/S}$ and $W_{S/T}$ are closed algebras of functions, and are equal.

Let $\Phi$ and  $\Psi$ be two-sided F{\o}lner sequences in $G.$  Let $f_1,f_2,f_3\in L^\infty(\mu),$ and define a sequence
\begin{align*}
A_n:=\frac{1}{|\Phi_n||\Psi_n|}\sum_{(g,h)\in \Phi_n\times \Psi_n} T_gf_1S_hf_2T_gS_hf_3.
\end{align*}
and auxiliary sequences
\begin{align*}
B_n^{(1)}&:= \frac{1}{|\Phi_n||\Psi_n|}\sum_{(g,h)\in \Phi_n\times \Psi_n} T_g(Pf_1)S_hf_2T_gS_hf_3\\
B_n^{(2)}&:= \frac{1}{|\Phi_n||\Psi_n|}\sum_{(g,h)\in \Phi_n\times \Psi_n} T_gf_1S_h(Qf_2)T_gS_hf_3\\
C_n&:= \frac{1}{|\Phi_n||\Psi_n|}\sum_{(g,h)\in \Phi_n\times \Psi_n} T_g(Pf_1)S_h(Qf_2)T_gS_hf_3.
\end{align*}

\begin{lemma}\label{red}  For $i=1,2,$ we have $\lim_{n\to \infty} A_n-B_n^{(i)}=0,$ and $\lim_{n\to \infty} A_n-C_n=0.$
\end{lemma}

\noindent\textit{Proof.}  We prove only that $\lim_{n\to \infty} A_n-B_n^{(1)}=0.$  The case $i=2$ is similar, and the last assertion follows from the first two.

Without loss of generality, suppose that $Pf_1=0.$  We need to show that $\lim_{n\to \infty} A_n=0,$ so we apply Lemma \ref{vdc'}.   For $g,h\in G$ let $u_{g,h}=T_gf_1S_hf_2T_gS_hf_3.$  Then for $g,j,j',h,k,k'\in G,$
\begin{align}\label{red1}
\langle u_{j'g,k'h}, u_{jg,kh} \rangle=\int T_g(T_{j'}f_1T_{j}\bar f_1)S_h(S_{k'}f_2S_k\bar f_2)T_gS_h(T_{j'}S_{k'}f_3T_jS_k\bar f_3)\, d\mu.
\end{align}
Applying $T_g^{-1}S_h^{-1}$ to the integrand in (\ref{red1}) yields
\begin{align}\label{red2}
\langle u_{j'g,k'h},u_{jg,kh}\rangle =\int S_h^{-1}(T_{j'}f_1T_j\bar f_1)T_g^{-1}(S_{k'}f_2S_k\bar f_2) (T_{j'}S_{k'}f_3T_jS_k\bar f_3)\, d\mu.
\end{align}
Averaging over $\Phi\times \Psi$  in the variables $g$ and $h$ we get
\begin{align}\label{red3}
\begin{split} \lim_{n\to\infty} \frac{1}{|\Phi_n\times \Psi_n|} &\sum_{(g,h)\in \Phi_n\times \Psi_n} \langle u_{j'g,k'h},u_{jg,kh}\rangle\\ &=\int \mathbb E(T_{j'}f_1T_j\bar f_1|\mathcal I_S)\mathbb E(S_{k'}f_2S_k\bar f_2|\mathcal I_T) T_{j'}S_{k'}f_3T_jS_k\bar f_3\, d\mu.
\end{split}
\end{align}
Applying the Cauchy-Schwartz inequality to the right-hand side of (\ref{red3}), and noting that $\|\mathbb E(S_{k'}f_2S_k\bar f_2|\mathcal I_T) T_{j'}S_{k'}f_3T_jS_k\bar f_3\|$ is bounded by some nonnegative constant $\alpha,$ we have the inequality
\begin{align}\label{red4}
\left| \lim_{n\to \infty} \frac{1}{|\Phi_n||\Psi_n|} \sum_{(g,h)\in \Phi_n\times \Psi_n}\langle u_{j'g,k'h},u_{jg,kh}\rangle\right| \leq \alpha \|\mathbb E(T_{j'}f_1T_j\bar f_1|\mathcal I_s)\|.
\end{align}
Averaging over $j,j'\in \Phi_m$ on both sides of (\ref{red4}) and applying Lemma (\ref{wm}), we get
\begin{align*}
\lim_{m\to \infty} \frac{1}{|\Phi_m|^2|\Psi_m|^2} \sum_{j,j'\in \Phi_m} \sum_{k,k'\in \Psi_m} \lim_{n\to \infty} \frac{1}{|\Phi_n||\Psi_n|}\sum_{(g,h)\in \Phi_n\times \Psi_n} \langle u_{j'g,k'h},u_{jg,kh}\rangle=0.
\end{align*}
By Lemma \ref{vdc'}, we have $\lim_{n\to \infty} A_n=0.$ $\square$
\medbreak

\subsection{Construction of an invariant measure on $X\times X\times X$.}   Lemma \ref{red} allows us to assume that $f_1\in W_{T/S}$ and $f_2\in W_{S/T}$ when we consider limits of the averages in (\ref{lim}). In this subsection, we show that this assumption allows us to approximate the summands $T_gf_1S_hf_2T_gS_h f_3$ in (\ref{lim}) by the linear combinations of continuous linear images of sequences the form $(T_g\times \text{id} \times T_g)(\text{id}\times S_h\times S_h)(\phi_1\otimes \phi_2\otimes \phi_3),$ where  $\phi_1\otimes \phi_2\otimes \phi_3$ is defined on $X\times X\times X$ with a $T\times \text{id} \times T$- and $\text{id}\times S\times S$-invariant Borel measure.  To be precise, there is a measure space $(Y,\mathcal D,\nu),$ measure preserving actions $T',S'$ of $G$ on $Y,$ and functions $\phi_1,\phi_2,\phi_3\in L^\infty(\nu),$ and a continuous linear map $M: L^2(\nu)\to L^2(\mu)$ so that the sequence $u_{g,h}:=T_gf_1S_hf_2T_gS_hf_3$ can be approximated by a linear combination of sequences of the form $v_{g,h}=M(T_g'\phi_1S_h'\phi_2T_g'S_h'\phi_3).$

For each $x\in X,$ let $\lambda_x=\mu_{\sigma(x)}\times\mu_{\tau(x)}\times \delta_x,$ where $\delta_x$ is the measure defined by $\delta_x(A)=1$ if and only if $x\in A.$  Let $\lambda=\int \lambda_x\, d\mu(x),$ so that $\lambda$ is a measure on $X\times X\times X.$  We may also define $\lambda$ by $\int f_1\otimes f_2\otimes f_3\, d\lambda =\int \mathbb E(f_1|\mathcal I_S)\mathbb E(f_2|\mathcal I_T) f_3\, d\mu.$  This second definition makes it clear that $\lambda$ is $T\times \text{id}\times T$-invariant and $\text{id}\times S\times S$-invariant.

We make the following standard observation about $\lambda$:
\begin{observation}\label{linear}  Let $H\in L^\infty(\lambda), f\in L^\infty(\lambda),$ and let $\{f_n\}_{n\in \mathbb N}$ be a uniformly bounded sequence of functions in $L^\infty(\lambda)$ such that $\lim_{n\to \infty} \|f_n-f\|_{L^2(\lambda)}=0.$  Let $F(x)=\int Hf\, d\lambda_x$ and for each $n,$ let $F_n(x)=\int Hf_n\, d\lambda_x,$ so that $F_n\in L^\infty(\mu).$  Then $\lim_{n\to \infty} \|F_n-F\|_{L^2(\mu)}=0.$
\end{observation}

We can now prove part (1) Theorem \ref{main}.
\noindent\textit{Proof of Theorem \ref{main}, part (1).}  Our goal is to show that the limit in (\ref{lim}) exists.  By Lemma \ref{red}, we may assume that $f_1\in W_{T/S}$ and $f_2\in W_{S/T}.$  Thus, $f_1$ may be approximated in $L^2(\mu)$ by linear combinations of functions of the form $\int H_1(x,z_1)\phi
_1(z_1) \, d\mu_{\sigma(x)}(z_1),$ where $H_1\in L^\infty(\mathbf X\underset{\mathbf I_S}\times \mathbf X)$ is $T\times T$-invariant and $\phi_1\in L^\infty(\mu),$ while  $f_2$ may be approximated in $L^2(\mu)$ by linear combinations of functions of the form $\int H_2(x,z_2) \phi_2(z_2)\, d\mu_{\tau(x)}(z_2),$ where $H_2\in L^\infty(\mathbf X\underset{\mathbf I_T}\times \mathbf X)$ is $S\times S$-invariant.  Since each $f_i$ is uniformly bounded, and each of the approximations to the $f_i$ are uniformly bounded, we may deduce the existence of the limit in (\ref{lim}) from the existence of limits of the form
\begin{align}\label{proof1}
\begin{split}\lim_{n\to \infty}\frac{1}{|\Phi_n||\Psi_n|} \sum_{(g,h)\in \Phi_n\times\Psi_n} &\left(\int H_1(T_gx,z_1)\phi_1(z_1)\, d\mu_{\sigma(T_gx)}(z_1) \right) \times \\ &\left(\int H_2(S_hx,z_2)\phi_2(z_2)\, d\mu_{\tau(S_hx)}(z_2)\right) f_3(T_gS_hx).
\end{split}
\end{align}
Write $u_{g,h}$ for the summand in (\ref{proof1}).  Computing as in equation (\ref{composition}), we can rewrite $u_{g,h}$ as the linear image of a function in $L^\infty(\lambda)$:
\begin{align}\label{rewrite}
\begin{split} &u_{g,h}=\int H_1(T_gx,z_1)\phi_1(z_1)\, d\mu_{\sigma(T_gx)}(z_1) \int H_2(S_hx,z_2)\phi_2(z_2)\, d\mu_{\tau(S_hx)}f_3(T_gS_hx)\\
&=\int H_1(x,z_1)\phi_1(T_gz_1)\, d\mu_{\sigma(x)}(z_1) \int H_2(x,z_2)\phi_2(S_hz_2)\, d\mu_{\tau(x)}f_3(T_gS_hx)\\
&= \int H_1(x,z_1)H_2(x,z_2) \phi_1(T_gz_1)\phi_2(S_hz_2)f_3(T_gS_hz_3)\, d\mu_{\sigma(x)}\times \mu_{\tau(x)} \times \delta_x(z_1,z_2,z_3).
\end{split}
\end{align}
Define a function $K\in L^\infty(\lambda)$ by $K(z_1,z_2,z_3)=H(x,z_1)H(x,z_2)$ if $\sigma(z_1)=\sigma(x)$ and $\tau(z_2)=\tau(x).$  From the last line of (\ref{rewrite}) we have
\begin{align}
u_{g,h}(x)=\int K\cdot (T\times \text{id}\times T)_g(\text{id} \times S\times S)_h \phi_1\otimes \phi_2 \otimes f_3\, d\lambda_x.
\end{align}
Averaging over $\Phi\times \Psi,$ Theorem \ref{ergodic} implies that
\begin{align*}
F:=\lim_{n\to\infty} \frac{1}{|\Phi_n||\Psi_n|} \sum_{(g,h)\in \Phi_n\times\Psi_n} (T\times \text{id}\times T)_g(\text{id} \times S\times S)_h \phi_1\otimes \phi_2\otimes f_3
\end{align*}
exists in $L^2(\lambda)$ and is in fact the orthogonal projection of $\phi_1\otimes \phi_2\otimes f_3$ on the space of $T \times \text{id}\times T$- and $\text{id} \times S\times S$-invariant functions.  By Observation \ref{linear} we conclude that the average of $u_{g,h}$ converges, and
\begin{align*}
\lim_{n\to \infty} \frac{1}{|\Phi_n||\Psi_n|} \sum_{(g,h)\in \Phi_n\times\Psi_n} u_{g,h}(x)= \int K\cdot F\, d\lambda_x \ \ \text{(in $L^2(\mu)$)}
\end{align*}
This concludes the proof of part (1) of Theorem \ref{main}. $\square$
\medbreak

In order to identify the limit in (\ref{lim}) and prove a combinatorial corollary, we first show that limits over products of F{\o}lner sequences can be evaluated by taking iterated limits in each F{\o}lner sequence separately.

\noindent\textit{Claim.}  Let $\{u_{g,h}\}_{g,h\in G}$ be a bounded sequence indexed by $G\times G$ such that for all two-sided F{\o}lner sequences $\Phi,\Psi,$ the limits
\begin{align*}
L&=\lim_{n\to\infty} \frac{1}{|\Phi_n||\Psi_n|}\sum_{(g,h)\in \Phi_n\times\Psi_n} u_{g,h}\\
L'&= \lim_{n\to \infty} \frac1{|\Phi_n|}\sum_{g\in \Phi_n} \lim_{m\to \infty} \frac{1}{|\Psi_m|}\sum_{h\in \Psi_m} u_{g,h}
\end{align*}
exist and are independent of the choice of F{\o}lner sequences $\Phi,\Psi.$  Then $L=L'.$
\medbreak

\noindent\textit{Proof.}  Let $\Phi, \Psi$ be two-sided F{\o}lner sequences in $G.$  For each $n,$ choose $m_n$ so that for all $g\in \Phi_n, \|\frac{1}{|\Psi_{m_n}|}\sum_{h\in \Psi_{m_n}} u_{g,h}-\lim_{m\to \infty} \frac{1}{|\Psi_m|}\sum_{h\in \Psi_m}u_{g,h}\|<\frac{1}{n}.$  Then
\begin{align*}
L-L'&=\lim_{n\to \infty} \frac{1}{|\Phi_n||\Psi_{m_n}|} \sum_{(g,h)\in \Phi_n\times \Psi_{m_n}} u_{g,h}-\lim_{n\to \infty}\frac{1}{|\Phi_n|} \sum_{g\in \Phi_n} \lim_{m\to\infty} \frac{1}{|\Psi_m|} \sum_{h\in \Psi_m} u_{g,h}\\
&= \lim_{n\to \infty} \frac{1}{|\Phi_n|} \sum_{g\in \Phi_n}\left(\frac{1}{|\Psi_{m_n}|}\sum_{h\in \Psi_{m_n}}u_{g,h} -\lim_{m\to \infty} \frac{1}{|\Psi_m|} \sum_{h\in \Psi_m} u_{g,h} \right)\\
&=0
\end{align*}
This proves the claim. $\square$
\medbreak

To prove the inequality claimed in Corollary \ref{large}, we will apply a corollary of the ergodic theorem for amenable groups.  Khintchine (\cite{K}) used the following result in the case where $G=\mathbb Z$ to prove a strengthening of the Poincar{\'e} recurrence theorem.  The argument we give is due to Hopf (\cite{H}).

\begin{cor}\label{Kh}  Let $\mathbf X,G,T$ and $\Phi$ be as in the statement of Theorem \ref{ergodic}, and let $f\in L^2(\mu)$ be non-negative.   Then
\begin{align*}
\lim_{n\to \infty} \frac{1}{|\Phi_n|} \sum_{g\in \Phi_n} \int fT_g \, d\mu\geq \left(\int f\, d\mu\right)^2.
\end{align*}
\end{cor}
\noindent\textit{Proof.}  By the ergodic theorem for amenable groups, we have
\begin{align*}
\lim_{n\to \infty} \sum_{g\in \Phi_n} \int fT_gf\, d\mu&=\int f\mathbb E(f|\mathcal I_T)\, d\mu\\
&=\int \mathbb E(f|\mathcal I_T)^2\, d\mu\\
&\geq \left(\int \mathbb E(f|\mathcal I_T)\, d\mu \right)^2\\
&= \left(\int f \, d\mu \right)^2,
\end{align*}
where the inequality is an application of Cauchy-Schwarz. $\square$

\noindent\textit{Proof of Corollary \ref{large}.}  Write $L$ for the limit in (\ref{large1}).  By Theorem \ref{main}, the limit $L$ exists, and by the preceding discussion, we can evaluate as follows: 
\begin{align*}
L&=\lim_{n\to \infty} \frac{1}{|\Phi_n|}\sum_{g\in \Phi_n} \lim_{m\to \infty} \frac{1}{|\Psi_m|} \sum_{h\in \Psi_m} \int fT_gfS_h(fT_gf)\, d\mu\\
&\geq \lim_{n\to \infty} \frac{1}{|\Phi_n|}\sum_{g\in \Phi_n} \left(\int fT_gf\, d\mu \right)^2\\
&\geq \left(\lim_{n\to \infty} \frac{1}{|\Phi_n|}\sum_{g\in \Phi_n}\int fT_gf\, d\mu  \right)^2\\
&\geq \left(\left( \int f\, d\mu \right)^2\right)^2,
\end{align*}
where we have applied either Corollary \ref{Kh} or Cauchy-Schwarz in every step.
$\square$
\medbreak
Let $A\subseteq X$ with $\mu(A)>0.$  We can apply Corollary \ref{large} with $f=1_A$ to conclude that intersections $(A\cap T_g^{-1}A)\cap S_h^{-1}(A\cap T_g^{-1}A)$ will have large measure for many pairs $(g,h).$  For this and later corollaries involving syndeticity, we need a general implication relating syndeticity and F{\o}lner sequences.

\begin{lemma}\label{tosynd}  Let $G$ be an amenable group, and $S\subseteq G\times G.$  Then $S$ is both left and right syndetic if and only if for all two-sided F{\o}lner sequences $\Phi,\Psi$ in $G,$ there exists $n,m\in \mathbb N$ such that $S\cap (\Phi_n\times \Psi_m)\neq \emptyset.$
\end{lemma}

\noindent\textit{Proof.}  Suppose that $S$ is not right syndetic.  Let $\Phi$ be a F{\o}lner sequence in $G.$  Then for all $n,m$ there exists $(x_n,y_m)\in G\times G\setminus \bigcup_{g\in \Phi_n\times \Phi_m} Sg^{-1},$ so $(x_n\Phi_n\times y_m\Phi_m)\cap S=\emptyset$ for all $n,m,$ which would contradict an assumption that $S$ meets $x_n\Phi_n\times y_m\Phi_m$ for some $n,m.$

Now suppose that $(\Phi_n\times \Psi_m)\cap S=\emptyset$ for all $n,m.$  Then $\Phi\times \Psi$ is a F{\o}lner sequence in $G\times G$ which does not meet $S,$ so $S$ is not syndetic.  $\square$

\begin{cor}\label{syndetic}
Let $(X,\mathcal B,\mu),T$ and $S$ be as in Theorem \ref{main}, and let $A\in \mathcal B.$  For all $\varepsilon>0$ the set
\begin{align*}
R_{\varepsilon}:=\{(g,h)\in G\times G:\mu((A\cap T_g^{-1}A)\cap S_h^{-1}(A\cap T_g^{-1}A))>\mu(A)^4-\varepsilon\}
\end{align*}
is both left- and right syndetic.
\end{cor}
\noindent\textit{Proof.}
Let $\varepsilon>0,$ and let $\Phi,\Psi$ be two-sided F{\o}lner sequences in $G.$  Take $f=1_A,$ so that $\mu((A\cap T_g^{-1}A)\cap S_h^{-1}(A\cap T_g^{-1}A))=\int fT_gfS_hfT_gS_hf \, d\mu.$ By Corollary \ref{large}, there exists $n$ so that
\begin{align*}
\frac{1}{|\Phi_n||\Psi_n|}\sum_{(g,h)\in \Phi_n\times \Psi_n}\mu(A\cap T_g^{-1}A\cap S_h^{-1}(A\cap T_g^{-1}A))>\mu(A)^4-\varepsilon.
\end{align*}
Since $\Phi$ and $\Psi$ are arbitrary two-sided F{\o}lner sequences, this shows that $R_{\varepsilon}$ is syndetic.
$\square$
\medbreak
\subsection{Identifying the limit; proof of part (2) of Theorem \ref{main}.}

Evaluating the iterated limit corresponding to (\ref{lim}) leads to an explicit description of the limit.  Recall that if $f\in L^\infty(\mu), \mathbb E(f|\mathcal I_S)(x)=\int f\, d\mu_{\sigma(x)}$ for $\mu$-almost every $x.$  For $f_i\in L^\infty(\mu)$ we have
\begin{align*}
\lim_{n\to \infty} \frac{1}{|\Phi_n||\Psi_n|} \sum_{(g,h)\in \Phi_n \times \Psi_n}& T_g f_1S_h(f_2 T_gf_3)(x)\\&= \lim_{n\to \infty} \frac{1}{|\Phi_n|} \sum_{g\in \Phi_n} \lim_{m\to \infty} \frac{1}{|\Psi_m|}\sum_{h\in \Psi_m}T_g f_1S_h(f_2 T_gf_3)(x)\\
&= \lim_{n\to \infty} \frac{1}{|\Phi_n|} \sum_{g\in \Phi_n} T_g f_1(x)\int f_2T_gf_3\, d\mu_{\sigma(x)}\\
&=\lim_{n\to \infty} \frac{1}{|\Phi_n|} \sum_{g\in \Phi_n} \int f_2(z)T_g\times T_gf_1\otimes f_3(x,z)\, d\mu_{\sigma(x)}(z)\\
&= \int f_2(z) H(x,z)\, d\mu_{\sigma(x)},
\end{align*}
where $H$ is the orthogonal projection of $f_1\otimes f_3$ on the space of $T\times T$-invariant functions in $L^2(\mu \underset{\mathcal I_S}\times\mathbf \mu).$  By symmetry, we find that the limit $L$ is also equal to $\int f_1(z)K(x,z)\, d\mu_{\tau(x)},$ where $K$ is the orthogonal projection of $f_2\otimes f_3$ on the space of $S\times S$-invariant functions in $L^2(\mu \underset{\mathcal I_T}\times \mu).$  Recalling the definition of $W_{T/S}$ and $W_{S/T},$ and noting that functions of the form $f\otimes g$ for bounded $f,g,$  span a dense subset of $L^2(\eta)$ whenever $\eta$ is a measure on $(X\times X, \mathcal B \otimes \mathcal B),$ we conclude that the limits in (\ref{lim}) span a dense subset of $W_{T/S}$ and of $W_{S/T}.$ It follows that $W_{T/S}=W_{S/T}.$

The fact that the limits in (\ref{lim}) span $W_{T/S}$ lets us conclude that limits will be constant almost everywhere exactly when the space $W_{T/S}$ consists solely of constant functions.  From the definition of $W_{T/S}$ and $W_{S/T},$ we see that this condition is equivalent to the condition that both $(X\times X, \mathcal B\otimes \mathcal B, \mu\times \mu, S\times S)$ and $(X\times X, \mathcal B\otimes \mathcal B, \mu\times \mu, T\times T)$ are ergodic systems - that is, that both $S\times S$ and $T\times T$ are weakly mixing.

\noindent\textit{Remarks.} (1) This description of the limit allows us to conclude that the factor determined by $W_{T/S}$ is the smallest characteristic factor for the scheme $(T_g,S_h,T_gS_h),$ in the sense that $W_{T/S}$ is a factor of every characteristic factor for the scheme $(T_g,S_h,T_gS_h).$ 
\medskip

(2)  The use of the measures $\lambda_x$ and $\lambda$ is essentially a ``diagonal measures" argument, as in the proof of Corollary 2.4 of \cite{F1}.  cf. \cite{A}.

\medskip

(3) The method of proof of Theorem \ref{main} may be used to prove the following similar result, which appears as Theorem 4.8 in \cite{BMZ}.

\begin{theo} Suppose that $G$ is a countable amenable group and $\Phi$ is a left F{\o}lner sequence for $G.$ Suppose that $(X,\mathcal B,\mu)$ is a probability space and that $T$ and $S$ are commuting measure preserving $G$-actions on $X.$  Then for any $\varphi,\psi\in L^2(X,\mathcal B,\mu),$
\begin{align*}
\lim_{n\to \infty} \frac{1}{|\Phi_n|} \sum_{g\in \Phi_n} \phi(T_gx)\psi(S_gT_gx)
\end{align*}
exists in $L^1(X,\mathcal B,\mu).$
\end{theo}

\section{Combinatorial Corollaries}\label{Combinatorial Corollaries}

We derive combinatorial corollaries analogous to those in sections 6 and 7 of \cite{BMZ}.

Let $\Omega=\{0,1\}^{G\times G}$ with the product topology, so that $\Omega$ is a compact metric space.  Regard elements of $\Omega$ as functions $1_E,$ where $E\subseteq G\times G.$  Define commuting $G$-actions $T$ and $S$ on $\Omega$ by $(T_g\xi)(g_1,g_2)=\xi(g_1g,g_2)$ and $(S_g\xi)(g_1,g_2)=(g_1,g_2g).$  

The following appears as Proposition 6.2 in \cite{BMZ}.

\begin{prop}\label{correspondence}  Suppose that $G$ is a countable amenable group and $\Phi$ is a left F{\o}lner sequence for $G\times G$.  Suppose $S\subseteq G\times G.$  Let $X=\overline{\{T_gS_h1_E:g,h\in G\}}$ be the orbit closure of $1_E$ in $\{0,1\}^{G\times G}.$  If
\begin{align*}
\bar d(E)=\limsup_{n\to \infty} \frac{|E\cap \Phi_n|}{|\Phi_n|}>0,
\end{align*}
then there exists a $\{T_g\}$- and $\{S_g\}$-invariant probability measure $\mu$ on $X$ such that
\begin{align*}
\mu(\{\eta\in X:\eta(e,e)=1\})=\bar d(E).
\end{align*}
\end{prop}
(The conclusion stated in \cite{BMZ} was that $\mu(\{\eta\in X:\eta(e,e)=1\})>0,$ but the above statement was in fact proved.)

The system constructed in Proposition \ref{correspondence} allows us to obtain Corollary \ref{combinatorial}, which is analogous to Theorem 6.1 in \cite{BMZ}.  To prove Corollary \ref{combinatorial}, let $E\subseteq G\times G,$ and take $(X,\mathcal B,\mu,T)$ to be the system constructed from $1_E$ in Proposition \ref{correspondence}.  The conclusion now follows from Corollary \ref{syndetic}.

\begin{cor}\label{combinatorial0} Let $E\subseteq G\times G,$ and let $\Phi$ be a left F{\o}lner sequence in $G\times G$ with $\limsup_{n\to \infty} \frac{|E\cap \Phi_n|}{|\Phi_n|}=\delta>0.$  Then for all $\varepsilon>0,$ the set
\begin{align*}
\left\{(g,h): \limsup_{n\to \infty} \frac{|E\cap E(g,1_G)\cap E(1_G,h)\cap E(g,h)\cap\Phi_n|}{|\Phi_n|}>\delta^4-\varepsilon\right\}
\end{align*}
is syndetic in $G\times G.$
\end{cor}

While most of the combinatorial corollaries of recurrence results in measure preserving dynamics pertain to configurations in sets of positive density in groups, Corollary 7.2 of \cite{BMZ} derives the following partition result from the recurrence statement in Theorem 5.2 of \cite{BMZ}.

\begin{theo} {\rm (\cite{BMZ}, Corollary 7.2)} Suppose that $G$ is a countable amenable group, $r\in \mathbb N,$ and that $G\times G\times G=\bigcup_{i=1}^r C_i.$ Then the set
\begin{align*}
\{g\in G: \text{ there exists } i, 1\leq i \leq r, &\text{ and } (a,b,c)\in G\times G\times G, \text{ such that }\\
&\{(a,b,c), (ga,b,c), (ga,gb,c), (ga,gb,gc)\} \subseteq C_i\}
\end{align*}
is both left and right syndetic in $G.$
\end{theo}
We derive an analogous result, concerning parallelepiped-like configurations in $G\times G\times G.$

\begin{cor}\label{partition}  Suppose that $G$ is a countable amenable group, $r\in \mathbb N,$ and that $G^3=\bigcup_{i=1}^r C_i.$  For $a=(a_1,a_2,a_3)\in G^3,$ let 
\begin{align*}
(g,h,k)\cdot (a_1,a_2,a_3)=\{&(a_1,a_2,a_3),(ga_1,a_2,a_3),(a_1,ha_2,a_3),(a_1,a_2,ka_3),\\
&(ga_1,ha_2,a_3),(ga_1,a_2,ka_3),
(a_1,ha_2,ka_3),(ga_1,ha_2,ka_3)\}.
\end{align*}  Then the set
\begin{align*}
\{(g,h,k)\in G\times G\times G: \text{ there exists } i, 1\leq i \leq r, \text{ and } a&\in G^3 \text{such that }\\
&(g,h,k) \cdot a \subseteq C_i\}
\end{align*}
is syndetic in $G\times G\times G.$  
\end{cor}

We will derive Corollary \ref{partition} from the following topological recurrence statement.

\begin{theo}\label{topological}  Let $R, S$ and $T$ be commuting actions of $G$ by homeomorphisms of a compact metric space $(X,d).$  Then for all $\varepsilon>0,$  the set
\begin{align*}J_\varepsilon=\{(g,h,k)\in & G\times G\times G: \text{ there exists } x\in X \text{ such that } \\ &\operatorname{diam}\{x,T_gx,S_hx,R_kx,T_gS_hx,T_gR_kx,S_hR_kx,T_gS_hR_kx\}<\varepsilon\}
\end{align*}
is both left and right syndetic in $G\times G\times G.$

\end{theo}

In the proof of Theorem \ref{topological} we will use the following fact: Given minimal topological system, $(X,d),$ with acting group $\Gamma,$ a point $x\in X$ and an open set $U\subset X,$ the set $\{\gamma \in \Gamma : \gamma x\in U\}$ is syndetic in $\Gamma.$  If $\Gamma$ is amenable, combining this fact with the ergodic theorem for amenable groups we see that every Borel probability measure on $X$ which is preserved by the action of $\Gamma$ assigns positive measure to every open set.

We remark that of the consequences of Corollary \ref{Kh} is that given a system $(X,\mathcal B,\mu,T),$ with $T$ a measure preserving action of an amenable group, and a set $A\in \mathcal B$ of positive measure, the set $\{g\in G: \mu(A\cap T_g^{-1}A)>0\}$ is both left- and right syndetic.

\noindent\textit{Proof.}  We may assume that the group generated by the individual elements of $R,S,$ and $T$ acts minimally on $X.$  

\noindent\textit{Claim.}  For all open $U\subseteq X,$ there exists a syndetic set $H\subseteq G\times G$ such that for all $(g,h)\in H,$ there exists $y\in U$ such that $\{y,T_gy,S_hy,T_gS_hy\}\subseteq U.$

To prove the claim, let $Y\subseteq X$ be a closed subset of $X$ minimal with respect to the group generated by the individual elements of $S$ and $T.$  Let $\nu$ be a measure on $Y$ invariant with respect to $S$ and $T.$  By minimality, we have, for all open $U$ with $U\cap Y\neq \emptyset,$ $\nu(U)>0.$  Fixing such a  $U$, we let $H$ be the set of $(g,h)$ such that $\nu(U\cap T_g^{-1}U\cap S_h^{-1}U\cap T_g^{-1}S_h^{-1} U)>0,$ so that $H$ is syndetic, by Corollary \ref{syndetic}.  Now suppose that $U$ is an arbitrary open subset of $X.$  By minimality of the action generated by the elements of $R,$ $S$ and $T,$ there exists $k\in G$ so that $Y\cap R_k^{-1}U\neq \emptyset,$ so by the above argument the set $H$ of $(g,h)$ such that $\nu(R_k^{-1}U\cap T_g^{-1}R_k^{-1}U\cap S_h^{-1}R_k^{-1}U \cap T_g^{-1} S_h^{-1}R_k^{-1}U)>0$ is syndetic, and so for all $(g,h)\in H, U\cap T_g^{-1}U\cap S_h^{-1}U\cap T_g^{-1}S_h^{-1} U\neq \emptyset.$  Thus, $H$ satisfies the conclusion of the claim.

Now let $\varepsilon>0.$  Let $\Phi,\Psi,$ and $\Theta$ be F{\o}lner sequences in $G.$  Let $U$ be an open set in $X$ of diameter less than $\varepsilon,$ and let $H\subseteq G\times G$ be a syndetic set such that for all $(g,h)\in H,$ there exists $y\in U$ with $\{y,T_gy,S_hy,T_gS_hy\}\subseteq U.$  Since $H$ is syndetic, there exists $n\in \mathbb N$ such that $H\cap (\Phi_n \times \Psi_n) \neq \emptyset.$  Fix such $(g,h)\in H\cap (\Phi_n \times \Psi_n)$ and $y\in U.$  Let $V$ be a neighborhood of $y$ such that $\{x,T_gx,S_hx,T_gS_hx\}\subseteq U$ for all $x\in V.$  By minimality of the action generated by $R,S,$ and $T,$ there is an $R,S,T$-invariant measure $\mu$ on $X$ with $\mu(V)>0.$  By Poincar\'e recurrence, $K:=\{k\in G: R_k^{-1}V\cap V\neq\emptyset \}$ is syndetic, so for all large enough $m,$ there exists $k\in K\cap \Theta_m.$  For such $k, x\in R_k^{-1}V\cap V,$ we have
$$
\{x,T_gx,S_hx,T_gS_hx,R_kx,T_gR_kx,S_hR_kx,T_gS_hR_kx\}\subseteq U.
$$
This shows that $J_\varepsilon$ meets $\Phi_n\times \Psi_n\times \Theta_m.$  Since $\Phi,\Psi,$ and $\Theta$ are arbitrary F{\o}lner sequences, Lemma \ref{tosynd} implies that $J_\varepsilon$ is both left and right syndetic. $\square$

To derive Corollary \ref{partition} from Theorem \ref{topological}, we apply a standard construction of a topological dynamical system from the partition $G\times G\times G=\bigcup_{i=1}^r C_i$ in the hypothesis of Corollary \ref{partition}.  Let $\Lambda=\{1,\dots, r\},$ and consider the shift space $X=(G\times G\times G)^\Lambda,$ with the product topology, together with the actions $R,S,T$ of $G$ defined by $(T_g \xi){x,y,z)}=\xi{(g^{-1}x,y,z)},(S_g \xi){x,y,z)}=\xi{(x,g^{-1}y,z)}, (R_g\xi){x,y,z)}=\xi{(x,y,g^{-1}z)}.$      Define a point $\xi_0$ in $X$ by $\xi(x,y,z)=i$ if $(x,y,z)\in C_i.$  Let $Y$ be the orbit closure of $x_0$ under the group generated by $R,S$ and $T.$  By Theorem \ref{topological}, the set of $(g,h,k)\in G^3$ such that there exists $\xi\in Y$ with $\operatorname{diam}\{\xi,T_g\xi, S_h\xi,  T_gS_h\xi, R_k \xi, T_gR_k\xi, S_hR_k\xi, T_gS_hR_k\xi \}<1$ is both left and right syndetic.  In particular, the set of such $(g,h,k)$ such that there exists $\xi$ such that every element of $\{\xi,T_g\xi,S_h\xi,T_gS_h\xi,R_k\xi, T_gR_k\xi,S_hR_k\xi, T_gS_hR_k\xi\}$ agrees at the $(e,e,e)$ coordianate is both left and right syndetic.  Taking inverses gives the conclusion of Corollary \ref{partition}.

\section{An example.}\label{example}  In this section we specialize to the case $G=\mathbb Z,$ and give an example of a class of systems where the characteristic factor for the scheme $(T_g,S_h,T_gS_h)$ is not an inverse limit of nilsystems.  This is in contrast to the case $S=T,$ where the characteristic factor (assuming ergodicity of $T$) is a compact group rotation, as Bergelson showed in \cite{B}.  We thank Alexander Leibman for suggesting this example.  (cf. the example presented in Section 3 of \cite{FK} .)

Let $\mathbf Y_0=(Y_0,\mathcal D_0,\mu_0,T_0),\mathbf Y_1=(Y_1,\mathcal D_1,\mu_1,S_1)$ be ergodic $\mathbb Z$-systems, and let $K$ be a compact metric group with Haar measure $m.$  Let $\tau:Y_0\to K, \sigma:Y_1\to K$ be measurable maps.  Define transformations $T$ and $S$ on $X:=Y_0\times Y_1\times K$ by $T(y_0,y_1,k)=(T_0y_0,y_1,\tau(y_0)k), S(y_0,y_1,k)=(y_0,S_1y_1,k\sigma(y_1)^{-1}).$  Then $T$ and $S$ commute, and preserve $\mu:=\mu_0\times \mu_1\times m.$  If $\tau$ is chosen so that $(y_0,k)\mapsto (T_0y_0,\tau(y_0)k)$ is ergodic, then the group generated by $S$ and $T$ acts ergodically on $X.$  One can verify that $W_{T/S}=L^2(\mu)$ in this situation, while such systems $(X,\mathcal B,\mu,T,S)$ are usually not isomorphic to inverse limits of nilsystems.  This is in sharp contrast to the characteristic factors for schemes involving powers of a single transformation, such as $(T^n,T^m,T^r,T^{n+m},T^{n+r},T^{m+r},T^{n+m+r}),$ which are shown in \cite{HK} to be inverse limits of nilsystems.

\end{document}